\documentclass{amsart}

\usepackage{epsfig,amsfonts,latexsym,amsmath,amssymb,amsthm,amstext,euscript}
\newtheorem{lemma}{Lemma}[section]
\newtheorem{theo}[lemma]{Theorem}
\newtheorem{coro}[lemma]{Corollary}
\newtheorem{rema}[lemma]{Remark}
\newtheorem{propos}[lemma]{Proposition}
\newtheorem{example}[lemma]{Example}
\def\bdem{\begin{proof}}
\def\edem{\end{proof}}
\def\bequ{\begin{equation}}
\def\eequ{\end{equation}}

\newcommand{\laba}{\label}

\newcommand{\rr}{\mbox{$    \rightarrow   $}}
\newcommand{\papa}{H^{\infty}}
\newcommand{\disc}{ {\Bbb D} }
\newcommand{\ov}{ \overline }
\newcommand{\noi}{\noindent}
\newcommand{\eiti}{e^{i\theta}}

\renewcommand{\Re}{\mbox{Re}\,}
\renewcommand{\Im}{\mbox{Im}\,}
\newcommand{\la}{\langle}
\newcommand{\ra}{\rangle}
\newcommand{\clospan}{ \ov{\mbox{span}} \, }
\def\span{\mathop{\rm span}\nolimits}

\newcommand{\Z}{\mathbb{Z}}
\newcommand{\C}{\mathbb{C}}
\newcommand{\G}{\mathcal{G}}

\newcommand{\R}{\mathbb{R}}
\newcommand{\hol}{\mathcal{H}}

\begin{document}

\title[{\sc Orbits of non-elliptic disc automorphisms}]
{Orbits of non-elliptic disc automorphisms on $H^p$}
\author{Eva A. Gallardo-Guti\'errez}
\address{Departamento de An\'alisis Matem\'atico, Facultad de Ciencias Matem\'aticas, Universidad Complutense de Madrid e IUMA,  Plaza de Ciencias 3, 28040 Madrid, Spain.}
\email{eva.gallardo@mat.ucm.es}
\author{Pamela Gorkin}
\address{Department of Mathematics, Bucknell University, PA 17837, Lewisburg, U.S.A.}
\email{pgorkin@bucknell.edu}
\author{Daniel Su\'arez}
\address{Departamento de Matem\'atica,  FCEyN, Universidad de Buenos Aires,  (1428) N\'{u}\~{n}ez, Cap. Fed.,  Argentina.}
\email{dsuarez@dm.uba.ar}



\date{March 2011, Revised version: October 2011}

% ----------------------------------------------------------------------

\begin{abstract}
Motivated by the \emph{Invariant Subspace Problem}, we describe
explicitly the closed subspace $H^2$ generated by the
limit points in the $H^2$ norm of the orbit of a thin
Blaschke product $B$ under composition operators $C_{\varphi}$
induced by non-elliptic automorphisms. This description exhibits a
surprising connection to model spaces.
Finally, we give a constructive characterization of the $C_{\varphi}$-eigenfunctions in $H^p$ for
$1\le p\le \infty$.
\end{abstract}

\maketitle

\section{Introduction and Preliminaries}

Given a bounded linear operator $T$ on a separable Banach space $\mathcal{B}$, the fact that an operator has invariant subspaces may not tell you much about the operator. On the other hand, knowing
that an operator has a large number of invariant closed subspaces,
and, in particular, a structured family, may make it possible to obtain a
lot of information about the action of the operator on
$\mathcal{B}$. In this context, it is helpful to focus on the behavior of the
cyclic subspaces generated by the elements of $\mathcal{B}$ under
$T$; in other words, focusing on the subspace
$$
K_{x}=\overline{\span \{T^{n} x:\; n\geq 0\}}^{\; \mathcal{B}},
$$
where $x\in \mathcal{B}$. Here $T^0$ denotes the identity operator
$I$  and $T^{n}=T\circ \stackrel{n\; \,
\mbox{\scriptsize{times}}}{\cdots}\circ \, T$.

It turns out that knowing the cyclic subspaces of
concrete operators on an infinite dimensional separable Hilbert
space $\mathcal{H}$ (operators that are \emph{universal}) would solve
the long-standing open question known as the Invariant Subspace Problem. Recall that an operator $U$ on
$\mathcal{H}$ is called universal, in the sense of Rota
\cite{Ro}, if for any bounded linear $T$ on $\mathcal{H}$, there
exists a complex constant $\lambda\neq 0$ and a closed invariant subspace
$\mathcal{M}$ of $U$ such that $U_{|_{\mathcal{M}}}$ is similar to $\lambda T$.

Note that  every bounded linear operator on an infinite
dimensional separable Hilbert space $\mathcal{H}$ would have a
non-trivial (closed) invariant subspace $\mathcal{M}$, that is,
$\mathcal{M}\neq \{0\}$ and  $\mathcal{M} \neq \mathcal{H}$, if
and only if the minimal invariant subspaces of a universal
operator $U$ on $\mathcal{H}$ are just one dimensional.
In the eighties, Nordgren, Rosenthal and Wintrobe \cite{NRW} proved
%, using Caradus Theorem,
that if $\varphi$ is a hyperbolic
automorphism of the unit disc and $\lambda$ is in the interior of
the spectrum of the composition operator $C_{\varphi}$ acting on
the classical Hardy space $H^2$, then
$C_{\varphi}-\lambda I$ is a universal operator on $H^2$
(see also \cite{NRW1}). Of course, the lattice of the closed invariant subspaces
of $C_{\varphi}-\lambda I$ coincides with that of $C_{\varphi}$.
Thus, it is important to study the closed invariant subspaces
of $C_{\varphi}$ in $H^2$ and, in particular, the cyclic
subspaces generated by $H^2$ functions.

Our discussion turns naturally to the factorization of $f\in H^p$ into its inner and outer factor.
The inner factor can be factored further, into a piece carrying all of %characterized by
its zeros (the Blaschke factor) and one with no zeros (the singular inner factor).
If the zero sequence of the Blaschke product is an interpolating sequence, the Blaschke product is said to be interpolating.
An important subclass of the interpolating Blaschke products is the set of thin Blaschke products; recall that
a Blaschke product $B$ with zeros $\{z_n\}_{n\ge 1}$ is said to be thin if $$\prod_{n\neq k} \rho(z_n,z_k) \rr 1$$ as
$k\rr\infty$, where $\rho$ denotes the pseudo-hyperbolic distance in the open unit disc $\mathbb{D}$.
When this holds, $\{ z_n \}$ is called  a thin sequence.
Thin Blaschke products have the closest behavior to finite Blaschke products that we can expect from infinite ones.

When $\varphi$ is a non-elliptic automorphism, in \cite{GG} the first two authors exhibited Blaschke products that are cyclic for $C_\varphi$ by
showing that the closed linear span of the limit points of their orbits is the whole space $H^2$. Clearly, this forces such Blaschke products to
be infinite, since the limit points of orbits of finite Blaschke products are constant functions of modulus $1$.
Here, we consider an arbitrary thin Blaschke product and characterize the closed linear span of limit points of its orbit, which is trivially invariant for $C_\varphi$.
Concretely, we prove the following (see Theorem \ref{limpace}):
Let $L_\varphi(B)$ denote the set of limit points, in the $H^2$ norm, of the orbit of a thin Blaschke product $B$ under the composition operator $C_\varphi$ where $\varphi$ is a non-elliptic automorphism. Then the
$H^2$-closure of the linear span of $L_\varphi (B)$
is either
$$H^2 \ \mbox{ or }\   (zb H^2)^\perp,
$$
where $b$ is a Blaschke product with simple (or no) zeros that satisfies $b\circ \varphi = \gamma\, b$ for some $\gamma\in\partial\disc$.
We then proceed to show that the same result holds for the Hardy spaces $H^p$ when $1\le p<\infty$.
We find the appearance of model spaces surprising, and we also see from this result that a natural question follows: What are
the eigenfunctions of $C_\varphi$? It is easy to see that $f$ is an eigenfunction if and only if the same is true for its Blaschke, singular and
outer factors.

In \cite{Cowen} Cowen studied eigenfunctions for composition operators.
Later, Matache  \cite{ma} characterized the singular inner
eigenfunctions of $C_\varphi$ in terms of the behavior of
pull-back measures (see also \cite{Mortini} for discrete singular inner functions).   Our approach provides separate characterizations
for Blaschke products, singular inner functions and outer
functions in $H^p$ that are eigenfunctions. The basic idea is the same in each of the three cases, but some of the technicalities are different.

The rest of the paper is organized as follows. Section
\ref{secBlaschke} is devoted to studying the orbit of thin
Blaschke products. In Section \ref{seceigenfunctions} we characterize the
eigenfunctions of composition operators induced by non-elliptic
disc automorphisms.
Finally, and for the sake of completeness, we end this
preliminary section by recalling some basic results and notation.

\subsection{Notation and basic results} Throughout this paper the open
unit disc of the complex plane will be denoted by $\mathbb {D}$
and $\partial \mathbb{D}$ will stand for its boundary. We will
denote  the space of holomorphic functions on $\mathbb{D}$ endowed
with the topology of uniform convergence on compacta by
$\mathcal{H}(\mathbb{D})$.

Recall that the Hardy space $H^p$, $1\leq p<\infty$,
consists of holomorphic functions $f$ on $\mathbb D$ for which the
norm
$$
\|f\|_p=\left ( \sup_{0\leq r<1} \int_{0}^{2\pi}
|f(re^{i\theta})|^p \, \frac{d\theta}{2\pi}\right )^{1/p}
$$
is finite. The space consisting of bounded analytic functions on
$\mathbb{D}$ endowed with the sup norm will be denoted by
$H^{\infty}$. A classical result due to Fatou states
that every Hardy function $f$ has non-tangential limit at $e^{i\theta}\in
\partial \mathbb{D}$, except possibly on a set Lebesgue measure
zero (see \cite{Du}, for instance). Throughout this work,
$f(e^{i\theta})$ will denote the non-tangential limit of $f$ at $e^{i\theta}$.

Recall that an automorphism $\varphi$ of $\mathbb{D}$ can be
expressed in the form
$$
\varphi(z)=e^{i\theta}\, \frac{p-z}{1-\overline{p}z} \qquad (z\in
\mathbb{D}),
$$
where $p\in \mathbb{D}$ and $-\pi<\theta\leq \pi$. Recall that
$\varphi$ is called  \emph{hyperbolic} if $|p|>\cos(\theta/2)$
(thus, $\varphi$ fixes two points on $\partial \mathbb{D}$);
\emph{parabolic} if $|p|=\cos(\theta/2)$ (so, $\varphi$ fixes just
one point, located on $\partial \mathbb{D}$) and \emph{elliptic} if $|p|<\cos(\theta/2)$
(therefore, $\varphi$ fixes two points in the Riemann sphere, one in $\disc$ and the other outside $\ov{\disc}$,
see \cite{Ah}, for example).

Throughout this paper, the involution that interchanges $0$ and $w$ will be denoted by
$$\varphi_w(z) = \frac{w-z}{1 - \overline{w}z },$$
where $z\in\disc$. For $w\in\disc$, we shall always denote this automorphism by $\varphi_w$. The pseudo-hyperbolic and hyperbolic
metrics for $z,w\in\disc$ are given, respectively, by
$$
\rho (z, w) = |\varphi_w(z)| \ \mbox{ and }\    \beta(z,w) =  \log  \frac {1 + \rho(z,w)}{1- \rho(z,w)},
$$
and we will denote by $D_\rho(z, r)$ and $D_\beta(z, R)$ the
respective closed balls of center $z$ and radius $r$, with $0\le r
<1$ and $R\ge 0$.

\begin{center}
\includegraphics[width=14cm]{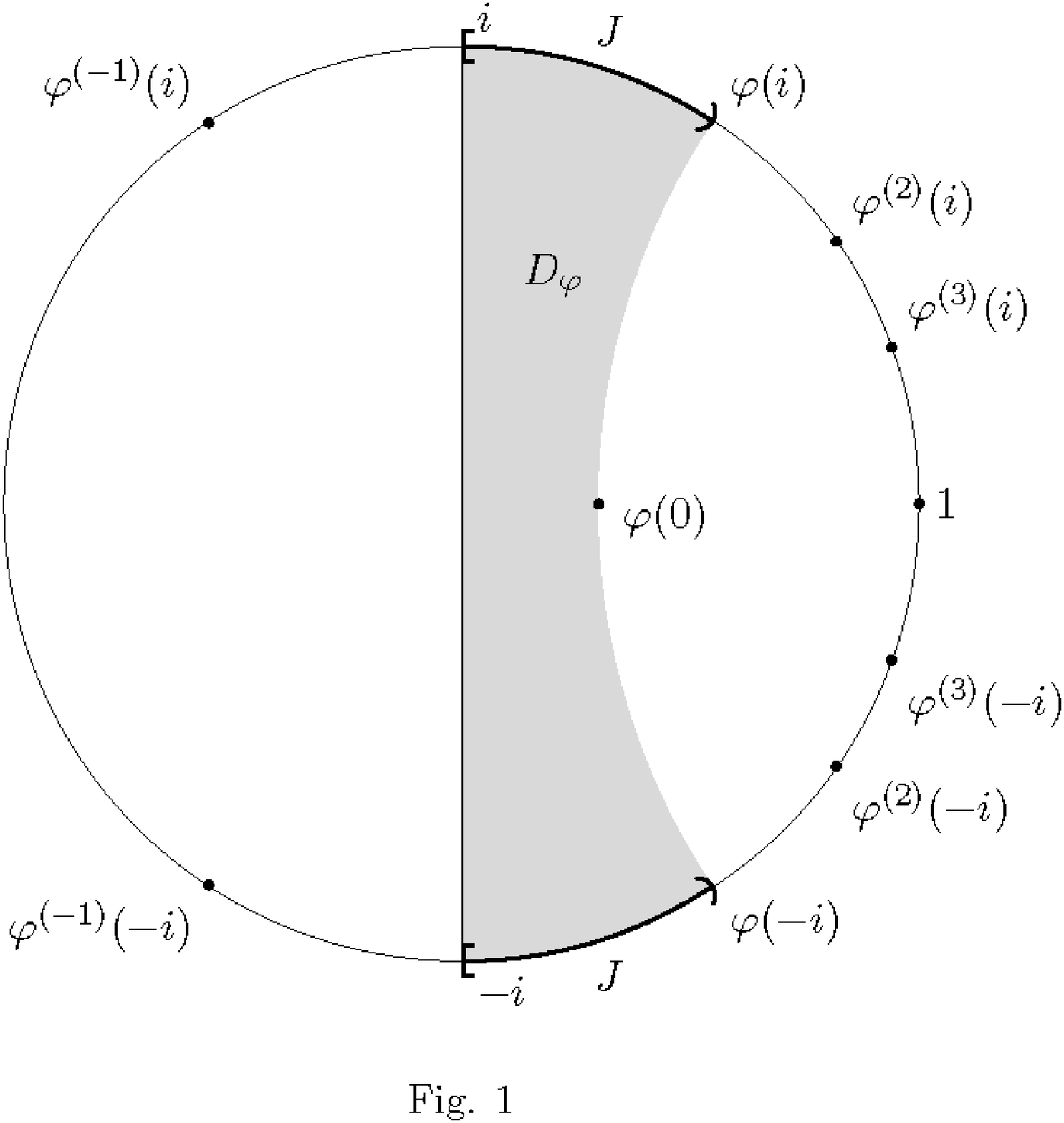}
\end{center}

\noi
If $\varphi$ is a non-elliptic automorphism and $n\in \Z$ (an
integer), we denote by $\varphi^{(n)}$  the $|n|$-th iterate of
$\varphi$ if $n>0$, of $\varphi^{-1}$ if $n<0$, and the identity
map if $n=0$. The action on $\disc$ of the group $\G_\varphi = \{
\varphi^{(n)}: n\in \Z \}$ leads naturally to the consideration of
the quotient space $\disc / \G_\varphi$, where we are identifying
points $z,w\in\disc$ such that $\varphi^{(n)}(z)=w$ for some
$n\in\Z$.

Since the class of $z\in\disc$ in $\disc / \G_\varphi$ is the
bilateral orbit $o_\varphi(z) =\{ \varphi^{(n)}(z): n\in \Z \}$,
we can represent  $\disc / \G_\varphi$ by any subset
$D_\varphi\subset \disc$ such that $o_\varphi(z) \cap D_\varphi$
is a singleton for every $z\in\disc$. Figures 1 and 2 show
reasonable choices of $D_\varphi$ when $\varphi$ is hyperbolic and
parabolic, respectively. We can transfer  the quotient topology of
$\disc/\G_\varphi$ to $D_\varphi$, so that the one-to-one
correspondence becomes a homeomorphism.

This identification allows us to think of the quotient map
$P : \disc \, \rr \, \disc/\G_\varphi \simeq D_\varphi$ as a map
onto $D_\varphi$; that is, $P(z)=  o_\varphi(z) \cap D_\varphi$,
defines a continuous map from $\disc$ onto $D_\varphi$, where
$D_\varphi$ has the quotient topology.
So, $D_\varphi$ can be identified with a subset of $\disc$, but endowed with the topology where a base of neighborhoods of a point
$z\in D_\varphi$ is given by $\{ w\in D_\varphi: \ \inf_{n\in\Z} |w-\varphi^{(n)}(z)| < \varepsilon \}$, for $\varepsilon >0$.
From Figures 1 and 2$\ $ it is not difficult to see that in both cases $D_\varphi$ is homeomorphic to a two-sided truncated cylinder
without the upper and lower boundaries.

\section{Orbits of thin Blaschke products} \label{secBlaschke}

We study the orbit of thin Blaschke products under composition
operators induced by non-elliptic automorphisms, characterizing
the closed set of its limit points (in the $H^2$ norm).

\begin{center}
\includegraphics[width=14cm]{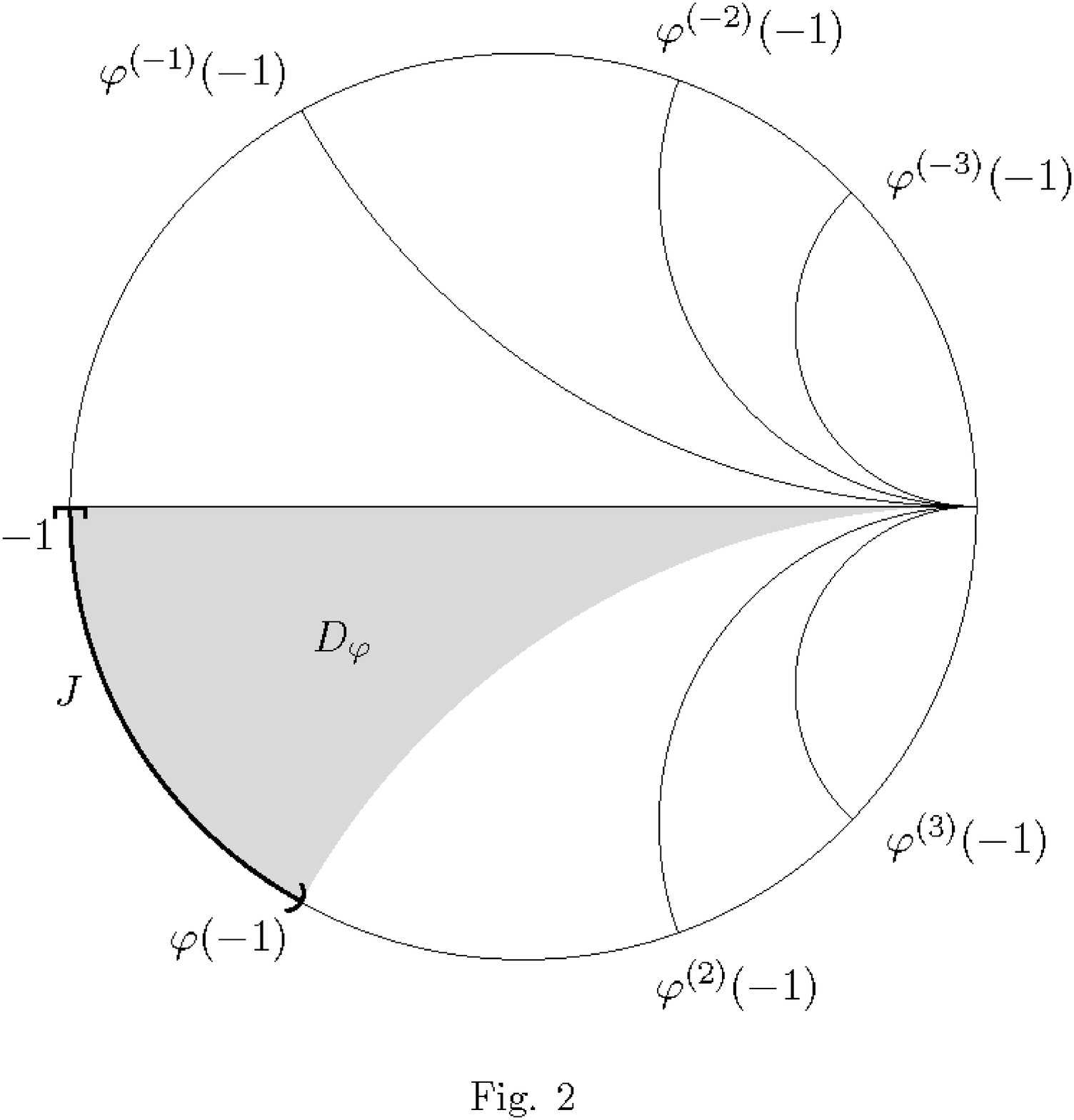}
\end{center}

\noi
If $\varphi$ is a self-map of $\disc$ and $f$ is analytic in
$\mathbb{D}$, that is, $f\in \hol(\disc)$, the orbit of $f$ is
defined by $O_\varphi(f)=\{ f\circ\varphi^{(n)} : n\ge 0\}$. When
$f$ is a bounded analytic function, that is, $f\in
H^\infty$, the orbit $O_\varphi(f)$ is precompact in the
compact-open topology of $\hol(\disc)$, so it makes sense to
define the limit set
$$
L_\varphi (f) = \{ g\in H^{\infty}: \
f\circ\varphi^{(n_k)} \rr g \,\mbox{ in } \hol(\disc) \mbox{ for
some subsequence $\{n_k\}$} \} .
$$
By the Corollary that follows Proposition 2 in \cite{BSh}, any sup-norm bounded sequence that converges in
$\hol(\disc)$ also converges weakly in $H^p$ for $1< p <\infty$.
Thus the points in $L_{\varphi}(f)$ belong to the
$H^p$-closure of the convex hull of $O_\varphi(f)$. As
we will show later, more can be said in case $f$ is a thin Blaschke product.
We proceed with a technical lemma that will be needed for our description of $L_{\varphi}(B)$ when $B$ is a thin Blaschke product.

\begin{lemma}\laba{sofax}
Let $\{z_n\}$ in $\disc$ be a sequence such that $|z_n| \rr 1$,  $0\leq r_k <1$ be any sequence and $0< \delta_k <1$ be a sequence that tends to $1$.
Then there is a subsequence $\{z_{n_k}\}$ of $\{z_n\}$ such that any sequence of points $\xi_k \in D_\rho(z_{n_k}, r_k)$ satisfies
$$
\prod_{j\neq k} \rho(\xi_j, \xi_k) > \delta_k .
$$
In particular, $\{\xi_k\}$ is thin.
\end{lemma}
\bdem Let $x_j \in (0,1)$ be a thin sequence such that $\prod_{j\neq k} \rho(x_j, x_k) > \delta_k$.
It will be enough to choose the points $z_{n_k}$ so that $\rho(D_\rho(z_{n_k}, r_k),
D_\rho(z_{n_j}, r_j)) \geq \rho(x_k,
x_j)\stackrel{\mbox{\scriptsize{def}}}{=} \alpha_{k,j}$ for all
$k\neq j$ or, equivalently,
$$
\beta(D_\beta(z_{n_k}, R_k), D_\beta(z_{n_j}, R_j)) \geq
\beta_{k,j} \ \mbox{ for } j<k,
$$
where
$$
R_k = \log \frac{1+r_k}{1-r_k}\ \mbox{ and }\ \beta_{k,j}=   \log
\frac{1+\alpha_{k,j}}{1-\alpha_{k,j}} .
$$
This is quite easy to do inductively. Once $z_{n_1}, \ldots,
z_{n_{k-1}}$ are chosen, simply take $z_{n_{k}}$ such that $\beta(
z_{n_{k}} , \{ z_{n_1}, \ldots, z_{n_{k-1}}  \} ) > R_k +
\sum_{j=1}^{k-1} (R_j+\beta_{k,j})$.
\edem

The preceding lemma resembles that of Wolff and Sundberg \cite[Lemma 5.4]{Su-Wo}, though we do not need to control the pseudo-hyperbolic distance from our given
sequence as they do. Instead, we pass to a subsequence to obtain the properties we need.
Our next result provides a description of all the sets in $\hol(\disc)$ of the form $L_\varphi(B)$, where $B$ is a thin Blaschke product.
The description contains some undetermined constants that will be irrelevant later when we take linear combinations
(in the proof of Theorem \ref{limpace}).

\begin{propos}\laba{iny}
Let $\varphi$ be a non-elliptic automorphism. If $B$ is a thin
Blaschke product, then there exists a nonempty set $V \subseteq
\partial \disc$ such that \bequ\laba{satur} L_\varphi(B) = \{
\lambda\varphi_w: \ w\in E , \,  \mbox{ for some } \lambda\in \partial\disc\}
\cup V, \eequ where $E\subset \disc$ is closed in $\disc$
with $\varphi(E)=E$. Conversely, given any such set $E$, there is
a thin product $B$ and a set $V$ such that $L_\varphi(B)$ is given by \eqref{satur}. Also, given $B$,
$$
E= \{ w\in\disc : \exists z_j\in Z(B) \mbox{ and integers } m_j \rr +\infty   \mbox{ such that  }
\varphi^{(-m_j)}(z_j) \rr w \} .
$$
\end{propos}

\bdem
Observe in advance that for $B$ thin, every $\hol(\disc)$-convergent subsequence of
$\{ B \circ \varphi^{ (n) } \}$ tends either to $\lambda \varphi_{w}$ or $\lambda$, for some
$\lambda\in \partial \disc$ and $w\in\disc$. This is because if $B \circ \varphi^{(n_{k})} \to f$ in
$\mathcal{H}(\mathbb{D})$, there are three possible situations: in the first one, $|f(0)| = 1$ and
consequently $B \circ \varphi^{(n_{k})} \to \lambda \in \partial{\mathbb{D}}$. If $f(0) = 0$, then the definition of thin
Blaschke product along with Schwarz's lemma  shows that there exists $\lambda \in \partial \mathbb{D}$ such that
$f(z) = \lambda z$ (see \cite[Proposition 2.3]{Hedenmalm}, for instance).
Finally, if $\xi = f(0) \ne 0$ and $|\xi| < 1$, then a computation shows that $\varphi_{\xi} \circ B$ is thin and
$\varphi_{\xi} \circ B \circ \varphi^{(n_{k})} \to \gamma z$ for some $\gamma \in \partial \mathbb{D}$. Thus,
$B \circ \varphi^{(n_{k})} \to \varphi_{\xi} (\gamma z) = \lambda \varphi_w(z)$, for some $\lambda\in \partial \disc$ and $w\in\disc$.

First we show that if $B$ is a thin product, there is some $\lambda\in\partial \disc$ that is a limit point of
$\{ B\circ\varphi^{(n)} \}$. Otherwise, the maximum modulus principle
implies that $\sup_{n\ge 0} |B(\varphi^{(n)}(0))| < \alpha$ for
some $\alpha <1$. Thus, there exists $R>0$ such that
$$
\varphi^{(n)}(0)\in \{ z: |B(z)| < \alpha \} \subset \bigcup_{ v\in Z(B) } D_\beta( v, R)
$$
for all $n\ge 0$. This is true because a thin product $B$ satisfies $|B(z)| \rr 1$ as $\beta(z, Z(B)) \rr \infty$
(see \cite[Ch.$\,$X, Lemma 1.4]{Ga}).
Since $\{ \varphi^{(n)}(0) \}$ is an interpolating sequence (see Section \ref{seceigenfunctions}),
the number of points contained in each of these balls, $D_\beta( v, R)$, must be bounded independently of the ball, say by $m$.
Consider only the zeros $v_k$ of $B$ such  that $D_\beta( v_k ,
R)$ contains some point $\varphi^{(n_k)}(0)$. Thus, at least one
of the points $\varphi^{(n_k)}(0), \varphi^{(n_k +1)}(0), \ldots ,
\varphi^{(n_k +m)}(0)$ must be contained in a different ball,
$D_\beta( v_j , R)$, with $j\neq k$. Hence,
\begin{eqnarray*}
\beta( v_k, v_j ) &\leq &\displaystyle \beta ( v_k , \varphi^{(n_k
)}(0) ) + \sum_{l=0}^{m-1} \beta(\varphi^{(n_k +l)}(0),
\varphi^{(n_k +l+1)}(0)) +R \\
&\leq & \displaystyle  R + \sum_{l=0}^{m-1}
\beta(\varphi^{(0)}(0), \varphi^{(1)}(0)) +R \\
&= & 2R +m  \beta(0, \varphi(0)),
\end{eqnarray*}
where the second inequality holds because $\beta$ is
$\varphi$-invariant, and this happens for every $k$. On the other
hand, since $\{ v_k \}$ is a thin sequence, $\beta ( v_k , \{ v_j
: j\neq k \} ) \rr \infty$ as $k\rr \infty$, a contradiction.

Now suppose that $w\in E= \big\{ w\in\disc: \lambda \varphi_w \in
L_\varphi(B) \mbox{ for some } \lambda\in\partial\disc \big\}$.
Then there is a sequence $\{n_k\}$ such that
$B\circ\varphi^{(n_k)} \rr \lambda \varphi_w$. So, for any
$m\in\Z$, $B\circ\varphi^{(n_k +m)} \rr \lambda (\varphi_w\circ
\varphi^{(m)} )$. Now, $\varphi_w\circ \varphi^{(m)}$  vanishes at
$\varphi^{(-m)} (w)$, so taking $m=-1, 1$ we see that
$\varphi^{(1)} (w), \, \varphi^{(-1)} (w)\in E$. Thus,
$\varphi(E)=E$. The fact that $E$ is closed in $\disc$
follows from a diagonal argument: indeed, if $w_s \in E$ is a sequence that tends to $w\in \disc$, then for each $s$ there exists a sequence $\{n_k(s)\}$ such that
$B\circ\varphi^{(n_k(s))} \rr \lambda_s \varphi_{w_s}$ for some $\lambda_s\in\partial\disc$, when $k\to \infty$.
We may assume that the sequence $\lambda_s$ converges to  some $\lambda_0 \in \partial \mathbb{D}$, and then we can extract a sequence from
$B\circ\varphi^{(n_k(s))}$ that tends to $\lambda_0 \varphi_w$, so $w\in E$.

To prove the converse, first consider the case in which $E=\emptyset$. Since
$\varphi^{(n)}(0) \rr \gamma \in \partial\disc$, the attractive fixed point of $\varphi$, for every
Blaschke product $B$ with zeros that do not accumulate at $\gamma$ we have
$B\circ \varphi^{(n)} \rr B(\gamma)\in\partial\disc$ in $\hol(\disc)$.

If $E\neq \emptyset$, choose a sequence $\{ \alpha_j \}$ that is
dense in $E$ and change it to $\{w_k\}_{k\ge 1}$, given by
$$
\alpha_1, \alpha_1, \alpha_2, \alpha_1, \alpha_2,
\alpha_3,\alpha_1, \ldots ,
$$
so that, as sets $\{w_k\}_{k\ge 1}= \{ \alpha_j \}_{j\ge 1}$, and
the set of limit points in $\disc$ of the sequence $\{w_k\}$ is
$E$. Write $\varphi^{(n)} = \lambda_n \varphi_{z_n}$, where
$|\lambda_n|=1$ and  $z_n\in \mathbb{D}$  (every automorphism can be written in this form), and observe that $|z_n|\to 1$.
Now use Lemma \ref{sofax} to choose a
subsequence $\{ z_{-n_k} \}$ of $\{ z_{-n}\}_{n\ge 0}$ such that
any sequence with one point in each $D_\rho(z_{-n_k}, \, |w_k|)$
is thin. Consider
$$
B(z) = \prod_{k\ge 1} \gamma_k(\varphi_{w_k} \circ
\varphi^{(-n_k)}) (z) = \prod_{k\ge 1} \gamma_k\varphi_{w_k}
(\lambda_{-n_k} \varphi_{z_{-n_k}} (z)),
$$
where $\gamma_k \in\partial\disc$ are chosen so that either each
factor is  positive at the origin or it is $z$. The zeros of $B$
are $\varphi_{z_{-n_k}}(\ov{\lambda}_{-n_k} w_k)\in
D_\rho(z_{-n_k}, \, |w_k|)$. If $w\in E$, there is a subsequence $w_{k_j} \rr w$, and
$$
B \circ \varphi^{(n_{k_j})}  = \varphi_{w_{k_j}} \gamma_{k_j}
\prod_{k\neq k_j} \gamma_k(\varphi_{w_k} \circ \varphi^{(-n_k +n_{k_j})}) .
$$
Since every $\hol(\disc)$-convergent subsequence of $\{ B \circ \varphi^{ (n_{k_j}) } \}$ tends
either to $\lambda \varphi_{z_\ast}$ or $\lambda$, for some $\lambda\in \partial \disc$ and $z_\ast\in\disc$,
and in our case $\varphi_{w_{k_j}}\rr \varphi_{w}$, it follows that
the convergent subsequences of $\{ B \circ \varphi^{(n_{k_j}) } \}$
tend to automorphisms of the form $\lambda\varphi_{w}$, for some $\lambda\in \partial \disc$.
On the other hand, if there is $\lambda\in\partial\disc$ such that $\lambda\varphi_w\in L_\varphi(B)$, we will show that $w\in E$.
For any integer $m$:
$$
Z(B\circ\varphi^{(m)}) = \varphi^{(-m)} (Z(B)) = \{ \varphi^{(-m+n_k)}(w_k): \ k\ge 1\} \subset E,
$$
where the last inclusion holds because $\{ w_k \}_{k\ge 1} \subset E$ and $\varphi(E)=E$. If $\lambda \varphi_w\in L_\varphi(B)$
for some $\lambda\in\partial\disc$ then $0=\inf_{m\ge 0} \rho(w, Z(B\circ\varphi^{(m)}))  \ge \rho(w,E)$, and consequently $w\in E$
because $E$ is closed.

For the last assertion of the proposition, suppose that $B$ is a thin product with  associated set $E$.
Then $w\in E$ if and only if  there exists a sequence $m_j\rr +\infty$ such that
$B\circ \varphi^{(m_j)} \rr \lambda  \varphi_w$ for some $\lambda \in \partial \disc$. This means that
$B( \varphi^{(m_j)}(w) ) \rr 0$, and since $B$ is an interpolating Blaschke product, this holds if and only if
$\rho(\varphi^{(m_j)}(w) , Z(B)) \rr 0$ (see \cite[p.$\,$395]{Ga}). That is, there is a sequence $\{z_j\}$ in $Z(B)$ such that
$\rho(\varphi^{(m_j)}(w) , z_j)=\rho(w ,\varphi^{(-m_j)}(z_j) ) \rr 0 $, or equivalently,
$w =\lim\varphi^{(-m_j)}(z_j)$.
\edem

\begin{lemma}\laba{limit}
Let $B$ be a thin Blaschke product. Then the $H^2$-limit points of $O_{\varphi}(B)$
are precisely the functions in $L_{\varphi}(B)$.
\end{lemma}

\bdem If $g$ is an $H^2$-limit point of $O_{\varphi}(B)$, there is a sequence
$B\circ \varphi^{(n_k)}\rr g$ in $H^2$.
It is clear that $g\in H^\infty$ and
bearing in mind that norm convergence in $H^2$ implies
uniform convergence on compacta, it follows that $g\in L_{\varphi}(B)$.

Conversely, if $g\in L_{\varphi}(B)$, there is a sequence $B\circ \varphi^{(n_k)}\rr g$ uniformly on compacta.
This, plus the fact that $\|B\circ \varphi^{(n_k)}\|_2=1$, implies that $B\circ\varphi^{(n_k)}\rr g$ weakly in
$H^2$.
By Proposition \ref{iny}, we have that either $g=\lambda \varphi_{w}$ or
$g=\lambda$, for some $\lambda\in
\partial \disc$ and $w\in
\mathbb{D}$. Hence,
\bequ\laba{trixie}
\|B\circ \varphi^{(n_k)}-g\|^2_2= 2-2 \Re \langle B\circ \varphi^{(n_k)}, g \rangle \to 2(1-\|g\|^2)= 0
\eequ
as $k\to \infty$, and the lemma follows. \edem

Now, we are in position to state the main result of this section.
If $N\subset H^2$ is a closed subspace, we denote  its orthogonal complement by $N^\perp$.

\begin{theo} \laba{limpace}
Let $\varphi$ be a non-elliptic automorphism and $B$ be a thin
product. Write $\mathcal{H}=\clospan L_\varphi (B)$ for the
closure in $H^2$ of the linear span of $L_\varphi (B)$.
Then
$$ \mathcal{H} = H^2 \ \mbox{ or }\   \mathcal{H}= (zb H^2)^\perp,
$$
where $b$ is a Blaschke product with simple (or no) zeros that
satisfies $b\circ \varphi = \gamma\, b$ for some
$\gamma\in\partial\disc$.
Conversely, if $\mathcal{H}$ is either of the above spaces, there
is a thin product $B$ such that $\clospan L_\varphi (B) =
\mathcal{H}$.
\end{theo}
\bdem Since Proposition \ref{iny} says that $L_\varphi (B)$
contains a non-null constant and for $w\in\disc\setminus \{0\}$,
\bequ\laba{angie} \ov{w}\varphi_w(z)= \ov{w} \frac{w-z}{1-\ov{w}z}
= 1- \frac{1-|w|^2}{1-\ov{w}z}, \eequ the same proposition tells
us that a subspace $\mathcal{H}\subset H^2$ has the form
$\mathcal{H}=\clospan L_\varphi (B)$ for $B$ thin if and only if
there is a relatively closed set $E\subset\disc$, with
$\varphi(E)=E$, such that
$$
\mathcal{H}= \clospan_{H^2} \left[ \{ 1 \} \cup \{
\frac{1}{1-\ov{w}z} : w\in E \} \cup \{ z : \mbox{ if } 0\in E \}
\right] .
$$
We recall that the function $K_w(z)= (1-\ov{w}z)^{-1}$ is the reproducing kernel in $H^2$ for $w\in\disc$, meaning that
$\la f, K_w \ra = f(w)$ for all $f\in H^2$.
So, a function $f$ is orthogonal to $\mathcal{H}$ if and only if $f(0)=0$, $f(w)= 0$ for all $w\in E$, and $f'(0)=0$ when $0\in E$.
That is,
$$
\mathcal{H}^\perp = \{ f\in H^2 : f\in zH^2, \ f\equiv 0 \ \mbox{ on } E \ \mbox{ and } f\in z^2 H^2 \ \mbox{ if } 0\in E \} .
$$
There are four possibilities: $E=\emptyset$, $E\neq\emptyset\, $ is not a Blaschke sequence, or $E$ is a Blaschke sequence,
in which case we distinguish temporarily between the cases $0\not\in E$ and $0\in E$.
If $E=\emptyset$, then $\mathcal{H}^\perp = zH^2$. Since the zeros of a non-null function in $H^2$ form a Blaschke sequence, when $E$
is not such sequence we have $\mathcal{H}^\perp = \{ 0 \}$. If $E$ is a Blaschke sequence, let $b$ and $b_0$ be Blaschke products that
satisfy $Z(b)= E$ and $Z(b_0) = E \setminus \{ 0\}$. If $0\not\in E$, we get
$\mathcal{H}^\perp = zH^2 \cap bH^2 = zbH^2$. On the other hand, if $0\in E$ we get
$\mathcal{H}^\perp = b_0H^2 \cap z^2 H^2 = z^2 b_0 H^2= zbH^2$.
Summing up,
\begin{enumerate}
\item[1.] $\mathcal{H}=H^2 \Leftrightarrow E$ is not a Blaschke sequence.
\item[2.] $\mathcal{H}=\C \Leftrightarrow E=\emptyset$.
\item[3.] $\mathcal{H}=(zb H^2)^\perp$, with $Z(b)=E$ $\Leftrightarrow E$ is a Blaschke
sequence, whether $0\in E$ or not.
\end{enumerate}
Finally, in the last case we have $Z(b\circ\varphi) = \varphi^{-1}(Z(b))= \varphi^{-1}(E)= E = Z(b)$. Therefore $b\circ \varphi$ and $b$ are Blaschke products
with the same zeros, which means that $b\circ \varphi= \gamma b$ for some unimodular constant $\gamma$.
\edem

By Theorem \ref{limpace}, the equality is attained in the inclusion
$\clospan L_\varphi(B) \subseteq \clospan O_\varphi(B)$
if and only if $\clospan L_\varphi(B)=H^2$ or
$\clospan L_\varphi(B)=(zb H^2)^\perp$, where $b$ is
a non-constant Blaschke product as in the theorem, and $B \in (zb H^2)^\perp$.
This means that $\ov{z}\ov{b}B = \ov{z}\ov{v}$, where $v\in H^2$ is inner,
or equivalently, that $B$ divides $b$.
Our next example shows that $\clospan L_\varphi(B)$ can be much
smaller than $\clospan O_\varphi(B)$.

\begin{example}
There exist thin Blaschke products $B$ such that
$$
\clospan O_\varphi (B)= H^2  \ \mbox{ and }\   \clospan L_\varphi
(B)= \mathbb{C}.
$$
\end{example}
\bdem Let $\{ w_k \}$ be a sequence in $D_\varphi$ such that $|w_k| \rr 1$ and whose set of non-tangential accumulation points has
positive measure. We show the existence of such a sequence as follows: take any maximal
sequence in $\disc$ that satisfies
$\beta (w_j , w_k) \ge \delta>0$ if $k \ne j$. Since the sequence is
hyperbolically separated, it does not accumulate on $\disc$.
We claim that it accumulates non-tangentially at every point of
$\partial\disc$.

Otherwise, for every angle $0<\alpha <\pi/2$, there is a point $\xi\in
\partial\disc$ such that some non-tangential sector with vertex $\xi$ and
half opening $\alpha$ does not contain any point $w_k$. If $\alpha$ is big
enough we have that the hyperbolic distance of $r\xi$ to the sequence is
larger than $\delta$ when $r<1$ is close to $1$. This contradicts the
maximality of the sequence. Hence, the intersection of this sequence with
$D_\varphi$ satisfies the desired condition.

As we did in the proof of Proposition \ref{iny}, we can use Lemma
\ref{sofax} to pick $n_k \rr +\infty$ so that the Blaschke product
$$
B(z) = \prod_{k\ge 1} \gamma_k (\varphi_{w_k} \circ
\varphi^{(-n_k)}) (z)
$$
with zeros $\varphi^{(n_k)} (w_k)$ is thin. Moreover, by Lemma \ref{sofax} we can choose $n_k \rr \infty$ so fast
that $B$ is as `thin' as we wish, meaning that if $j \ge 1$ and %%we denote
$B_j := \prod_{k\neq j}  \gamma_k (\varphi_{w_k} \circ \varphi^{(-n_k)})$,
then $|B_j(\varphi^{(n_j)} (w_j))|$ is so close to $1$ as we predefine
(by choosing $\delta_j \to 1$ fast enough in Lemma \ref{sofax}).
Furthermore, since for every $h\in H^\infty$, with $\|h\|_\infty \le 1$, and $w\in \disc$, the Schwarz-Pick inequality easily yields
$\rho(|h(w)|, |w|) \le |h(0)|$ (see \cite[Ch.$\,$I, Cor.$\,$1.3]{Ga}), taking $h= B_j\circ\varphi^{(n_j)}$ and $w=w_j$,
we can ensure that the right hand side of
$$\rho(|B_j \circ\varphi^{(n_j)}(w_j)|, |w_j|) \le |B_j \circ\varphi^{(n_j)}(0)|$$
is as close to 1 as desired.
In particular, we can impose the condition
\bequ\laba{uvan} \frac{(1-|B_j \circ\varphi^{(n_j)}(0)|^2  )^{ \frac{1}{2}} }{1-|w_j|^2} \,\rr\, 0 .
\eequ
We will check that $\clospan O_\varphi (B)= H^2$.
Suppose that $f\in H^2$ is orthogonal to $O_\varphi(B)$. By Proposition \ref{iny}, $f\perp \mathbb{C}$.
So, using \eqref{angie}, $f\perp (B\circ\varphi^{(n_k)})$, $f\perp \mathbb{C}$ and
$$B\circ\varphi^{(n_k)} = (B_{k}\circ\varphi^{(n_k)}) \gamma_k \varphi_{w_k}, \ (k \ge 1),$$
successively in the following chain of equalities, we get
\begin{eqnarray*}
\ov{{\gamma}_k B_{k}(\varphi^{(n_k)}(0))} \, f(w_k) & = &
\frac{1}{1-|w_k|^2} \left\la  f \, , \, \left(1 - \ov{w}_k\varphi_{w_k}\right) \gamma_k (B_{k}\circ\varphi^{(n_k)})(0) \right\ra  \\
& = &
\frac{1}{1-|w_k|^2}\left\la f \, , \, \ov{w}_k(B\circ\varphi^{(n_k)}) +   \left(1 - \ov{w}_k\varphi_{w_k}\right) \gamma_k (B_{k}\circ\varphi^{(n_k)})(0) \right\ra\\
&= & \left\la f \, , \, \frac{ \ov{w}_k}{1-|w_k|^2} (B\circ\varphi^{(n_k)}) - \frac{\gamma_k \ov{w}_k}{1-|w_k|^2} \varphi_{w_k} (B_{k}\circ\varphi^{(n_k)})(0) \right\ra \\
&=& \frac{w_k}{1-|w_k|^2} \, \la f \, , \,  \gamma_k \varphi_{w_k}  [(B_{k}\circ\varphi^{(n_k)})  -  (B_{k}\circ\varphi^{(n_k)})(0)] \ra .
\end{eqnarray*}
Consequently
\begin{eqnarray*}
|B_{k}(\varphi^{(n_k)}(0))| \, |f(w_k)| &\leq &
\frac{|w_k|}{1-|w_k|^2} \| f\|_2 \, \|(B_{k}\circ\varphi^{(n_k)})  -  (B_{k}\circ\varphi^{(n_k)})(0)\|_2 \\
&= & \frac{|w_k|}{1-|w_k|^2} \| f\|_2 \, (1-|(B_{k}
\circ\varphi^{(n_k)})(0)|^2  )^{ \frac{1}{2}} \rr 0
\end{eqnarray*}
by \eqref{uvan}; that is, $f(w_k) \rr 0$. Since the sequence $\{
w_k\}$ accumulates non-tangentially on a set of positive measure in
$\partial \disc$,  the non-tangential limits of $f$ must vanish on a
set of positive measure, implying that $f=0$. This proves our claim.

The equality $\clospan L_\varphi (B)= \mathbb{C}$ will follow from
the last assertion of Proposition \ref{iny} if we show that the set
$E$ associated with $B$ is empty. So, suppose that $w\in E$. Since
the zeros of $B$ are $\varphi^{(n_k)} (w_k)$, the proposition says
that there are integers $m_j, k_j \rr +\infty$ such that
$\varphi^{(-m_j)} ( \varphi^{ (n_{k_j}) } (w_{k_j}) ) \rr w$.
Applying the quotient map  $P: \disc \rr D_\varphi$, we obtain
$$
\lim w_{k_j} = \lim P \big( \varphi^{ (-m_j+n_{k_j}) }
(w_{k_j})\big)  = P(w) \in D_\varphi  ,
$$
where the limit is taken in the $D_\varphi$ topology. This is not possible, because $|w_{k_j}| \rr 1$.
\edem

\subsection{Generalization to Hardy spaces}

The proof of Theorem \ref{limpace} is more natural and transparent in the context of the Hilbert space $H^2$, but with very minor modifications
we can see that both Lemma \ref{limit} and Theorem \ref{limpace} are valid for $H^p$ when $1\le p <\infty$. We recall that the composition operator
$C_\varphi$ is bounded on $H^p$ for $p\ge 1$ (see \cite[p.$\,$121]{Co-Mac}).
To see that the lemma holds simply replace \eqref{trixie} by
$$
\|B\circ \varphi^{(n_k)}-g\|^p_p \le   \|B\circ \varphi^{(n_k)}-g\|^{p-1}_{2(p-1)}  \   \|B\circ \varphi^{(n_k)}-g\|_2
\le 2^{p-1}\   \|B\circ \varphi^{(n_k)}-g\|_2 ,
$$
which  by  \eqref{trixie} tends to $0$ as $k\to \infty$, where the first inequality is
the Cauchy-Schwarz inequality and the second holds because $|B\circ \varphi^{(n_k)}-g| \le 2$.

Let $1<p<\infty$ and $N\subset H^p$ be a subspace. The annihilator of $N$ is
$$N^\perp = \{ f \in H^q : \la h, f \ra =0 \mbox{ for all } h\in N \},$$
where  $1/p +1/q =1$.
Mimicking the proof of Theorem \ref{limpace}
we see that if $\varphi$ is a non-elliptic automorphism and $B$
is a thin product, a space $\mathcal{H} \subset H^p$ has the
form $\mathcal{H}=\clospan_{H^p} L_\varphi (B)$ if and only if it is closed and
$\mathcal{H}^\perp = \{ 0 \}\subset H^q\,$ or $\mathcal{H}^\perp = zb H^q$,
where $b$ is a Blaschke product as in Theorem \ref{limpace}.  This is the same as saying that
$$
\mathcal{H}= H^p \,  \ \mbox{ or }\  \,  \mathcal{H} = (zb H^q)^\perp = H^p \cap b \ov{H}^p,
$$
where the bar means complex conjugation.
This is the $H^p$ version of  Theorem \ref{limpace} for $1<p<\infty$.
The result for $H^1$ requires a different argument. Since the first of the following spaces
$$\span  L_\varphi (B) \subset  \clospan_{H^2} L_\varphi (B)   \subset       \clospan_{H^1} L_\varphi (B):= \mathcal{H}$$
is dense in the latter, so is the one in middle, which by Theorem \ref{limpace}, is either $H^2$ or $H^2 \cap b \ov{H}^2$,
with both cases ocurring.
In the first case $\mathcal{H}=H^1$ and in the second, $\mathcal{H}$ is the $H^1$-norm closure of $H^2 \cap b \ov{H}^2$,
which by Lemma 5.8.14 of \cite{Ci-Ro} is $H^1 \cap b \ov{H}^1$. This gives Theorem \ref{limpace} for $H^1$.

\section{A constructive characterization of eigenvectors} \label{seceigenfunctions}

While our work in this section is related to those in \cite[Section 4]{Cowen}, our approach to the problem is quite different.
The function $h(z) = i\frac{1+z}{1-z}$ maps $\disc$ onto
$\C_+ =\{v\in \C : \, \Im v >0 \}$, with $h^{-1}(v) =\frac{v-i}{v+i}$.
If $\varphi$ is a hyperbolic automorphism of $\disc$, by conjugating
$C_\varphi$ with an invertible operator, we can assume that its
fixed points are $-1$ and $1$, where $1$ is the attractive point.
Thus, $\varphi = h^{-1}\circ \tilde{\varphi} \circ h$, where
$\tilde{\varphi}: \C_+ \rr \; \C_+$ is $\tilde{\varphi}(w) = \alpha w$, with $\alpha >1$ (see \cite{Ah}, for instance).
A straightforward calculation shows that
\bequ\laba{mod-h}
\varphi(z) = - \varphi_{ \frac{1-\alpha}{1+\alpha} }(z) \ \mbox{ and  }\ \varphi^{(n)}(z) =
- \varphi_{ \frac{1-\alpha^n}{1+\alpha^n} }(z) \ \ \forall n\in\Z.
\eequ

\noi The same argument with $\varphi: \disc \rr \disc$ parabolic, where we now assume that its fixed point is $1$,
and therefore $\tilde{\varphi} (w)= w+t$, with $t\in \R\setminus \{0\}$, shows that
\bequ\laba{mod-p}
\varphi(z) = \left(\frac{t-2i}{t+2i}\right) \varphi_{ \frac{t}{t-2i} }(z) \ \mbox{ and  }\ \varphi^{(n)}(z) =
\left(\frac{nt-2i}{nt+2i}\right) \varphi_{ \frac{nt}{nt-2i} }(z) \ \ \forall n\in\Z .
\eequ
We can further assume that $t>0$, since otherwise the treatment is analogous.

From the above expressions for $\tilde{\varphi}$, it is easy to see in both cases that if $w_0\in  \C_+$ then  $\{ \tilde{\varphi}^{(n)}(w_0): \ n\in \Z \}$ is an interpolating sequence for
$\papa(\C_+)$ (see \cite[Ch.$\,$VII]{Ga}), so  $\{ \varphi^{(n)}(z_0): \ n\in \Z \}$ is an interpolating sequence for $\papa(\disc)$ for any fixed $z_0\in \disc$.

\subsection{Blaschke product eigenvectors}

\noi It is clear that a Blaschke product $b$ satisfies $b\circ
\varphi = \gamma \, b$ for some $\gamma\in\partial\disc$ if and
only if $b$ and $b\circ\varphi$ have the same zeros. This means
that for every zero $w$ of $b$, $\{ \varphi^{(n)}(w) : \,n\in\Z
\}$ are zeros of $b$, each one with the same multiplicity as $w$.
Hence,
$b$ is solely determined by the sequence of its zeros in $D_\varphi$.
So, a characterization of the Blaschke products that are
eigenvectors is tantamount to a characterization of the Blaschke
sequences $\{z_k\}_{k\ge 1}$ in $D_\varphi$ such that $\{
\varphi^{(n)}(z_k) : \,k\ge1, \,n\in\Z \}$ is also a Blaschke
sequence. We have to distinguish between the hyperbolic and the parabolic case.

\begin{theo}\laba{cabla}
Let  $\{z_k\}_{k\ge 1}$ be Blaschke sequence in $D_\varphi$. Then
$\{ \varphi^{(n)}(z_k) : \,k\ge1, \,n\in\Z \}$ is a Blaschke
sequence
\begin{enumerate}
\item[ (1)] always when $\varphi$ is hyperbolic.
\item[ (2)] if and only if $\{ z_k\}$ stays outside of some horocycle
tangent to $\partial \disc$ at the fixed point of $\varphi$, when
$\varphi$ is parabolic.
\end{enumerate}
\end{theo}
\bdem Geometrically, it will be more convenient to look at things
in the upper half-plane $\mathbb{C}_+ = \{ z\in \mathbb{C}: \Im z >0\}$. Hence, the Blaschke condition for the sequence
$z_k=x_k+iy_k$ is $\sum \frac{y_k}{1+|z_k|^2} < \infty$. To prove
(1) we can assume that $\varphi$ is a hyperbolic automorphism that
fixes $0$ and $\infty$, and therefore has the form  $\varphi(w)=
\alpha w$ with $1 \neq \alpha>0$. We can also assume that $\alpha
>1$, since the proof is the same in both cases. In this case, we
can take $D_\varphi = \{ z\in \mathbb{C}_+ : 1 \le |z| < \alpha
\}$, and we must show that if $\{z_k\}_{k\ge 1}$ in $D_\varphi$ is
a Blaschke sequence then so is $\{ \alpha^n z_k : k\geq 1, \, n\in
\mathbb{Z} \}$. Hence,
\begin{eqnarray*}
\sum_{k\ge 1, \, n\in \mathbb{Z}} \frac{\alpha^n y_k}{1+|\alpha^n z_k|^2} &\le &
\displaystyle  \sum_{n\ge 0}\sum_{k\ge 1}
\frac{\alpha^n y_k}{\alpha^{2n}  \, |z_k|^2} +  \sum_{n< 0}\sum_{k\ge 1} \alpha^n y_k \\
&\le & \displaystyle  \left[ \sum_{n\ge 0}  \frac{1}{\alpha^{n}} + \sum_{n< 0} \alpha^n \right]  \sum_{k\ge 0}  y_k
\hspace{22mm} \mbox{ (since $|z_k|\ge 1$)} \\
&\le & \displaystyle  \left[ \frac{\alpha+1}{\alpha-1}  \right]
\sum_{k\ge 0} \frac{ y_k(1+\alpha^2) }{1+|z_k|^2}  < \infty
\hspace{17,5mm} \mbox{ (since $|z_k|< \alpha$)} .
\end{eqnarray*}
To prove (2) we can assume that $\varphi$ is a parabolic
automorphism that fixes $\infty$, and therefore has the form
$\varphi(w)= w+t$, where $t\in\mathbb{R}$, and we can also assume
that $t>0$, since the proof is the same in both cases.

In the case at hand, we take $D_\varphi = \{ z\in \mathbb{C}_+ :
0\leq \Re z < t \}$, and show that given a Blaschke sequence
$\{z_k\}_{k\ge 1}$ in $D_\varphi$, $\{  z_k + nt : k\geq 1, \,
n\in \mathbb{Z} \}$ is a Blaschke sequence if and only if the
sequence $\{ y_k \}$ is bounded.

First suppose that there is some constant $C>0$ such that $y_k \le
C$ for all $k\ge 1$. Since $|x_k+nt| \ge |n|t- x_k \ge (|n|-1)t$,
we have that $(x_k+nt)^2 \ge (|n|-1)^2 t^2$ for all $n\neq 0$.
Consequently, if $n\in \mathbb{Z}\setminus \{ 0 \}$,
$$
\frac{1+|z_k|^2}{1+|z_k+nt|^2} \le    \frac{1+t^2+C^2}{1+(|n|-1)^2
t^2},
$$
leading to
\begin{eqnarray*}
\sum_{k\ge 1, \, n\in \mathbb{Z}\setminus \{ 0 \}}
\frac{y_k}{1+|z_k+nt|^2}  & = & \displaystyle
\sum_{k\ge 1, \, n\in \mathbb{Z}\setminus \{ 0 \}} \frac{y_k}{1+|z_k|^2} \, \frac{1+|z_k|^2}{1+|z_k+nt|^2} \\
& \le & \displaystyle   \sum_{k\ge 1} \frac{y_k}{1+|z_k|^2}
\sum_{n\in \mathbb{Z}\setminus \{ 0 \}}
\frac{1+t^2+C^2}{1+(|n|-1)^2 t^2},
\end{eqnarray*}
which is convergent. Obviously, the sum is also bounded for $n=0$.

If  the sequence  $\{ y_k \}$ is not bounded, fix an arbitrary
$y_k \geq t+1$ and consider all the values of $n\in \mathbb{Z}$
such that $-\frac{y_k}{t} \leq -\left[\frac{y_k}{t}\right]\leq n
\leq \left[\frac{y_k}{t}\right] \leq\frac{y_k}{t}$, where $[a]$
denotes the largest integer $\leq a$. Since $0\leq x_k <t$, we
have
$$
-y_k \le x_k-y_k \le x_k+nt \le x_k+y_k \le t+y_k \le 2y_k,
$$
and consequently $(x_k+nt)^2 \le (2y_k)^2$. Thus,
\begin{eqnarray*}
\sum_{-\left[\frac{y_k}{t}\right]\leq n \leq
\left[\frac{y_k}{t}\right]}  \frac{y_k}{1+|z_k+nt|^2} &\geq &
\frac{1}{y_k} \sum_{-\left[\frac{y_k}{t}\right]\leq n \leq
\left[\frac{y_k}{t}\right]} \frac{y^2_k}{1+ 5 y_k^2}\\
&\stackrel{\mbox{\scriptsize{($y_k \ge 1$) }}}{\ge}& \displaystyle
\frac{1}{y_k} \sum_{-\left[\frac{y_k}{t}\right]\leq n \leq
\left[\frac{y_k}{t}\right]} \frac{1}{6} \\
&=& \displaystyle  \frac{1}{6y_k} \left(
2\left[\frac{y_k}{t}\right] + 1 \right) \\
& \geq & \displaystyle  \frac{1}{6y_k} \left( \frac{y_k}{t}
\right) = \frac{1}{6t}.
\end{eqnarray*}
Consequently,
$$
\sum_{k\ge 1, \, n\in \mathbb{Z}} \frac{y_k}{1+|z_k+nt|^2} \ge
\sum_{y_k \ge t+1} \, \sum_{-\left[\frac{y_k}{t}\right]\leq n \leq
\left[\frac{y_k}{t}\right]} \frac{y_k}{1+|z_k+nt|^2} \ge \sum_{y_k
\ge t+1}  \frac{1}{6t},
$$
which is infinite if there are infinitely many $y_k \ge t+1$.
\edem

\subsection{Outer eigenfunctions}

\noi If $\varphi$ is given by \eqref{mod-h} (the hyperbolic case)
consider the circular intervals $J= (\varphi(i),i]\,\cup\, [-i,\varphi(-i))$ of
$\partial \disc$, and if $\varphi$ is given by
\eqref{mod-p} (the parabolic case) consider $J=[-1, \varphi(-1))$
(see Figures 1 and 2).

\begin{lemma}\laba{deriv}
If $\varphi$ is given by \eqref{mod-h} then \bequ\laba{des-h}
\frac{1}{4} \alpha^{-|n|} \le |\varphi^{(n)'}(w)| \le
(\alpha+1)^2 \alpha^{-|n|} \, \ \ \forall w\in J, \, \forall
n\in\Z. \eequ If $\varphi$ is given by \eqref{mod-p} then there is
a constant $c(t) >0$ such that \bequ\laba{des-p}
 \frac{1}{n^2 t^2+4} \leq  |\varphi^{(n)'}(w)|
\le \frac{c(t)}{n^2 t^2+4} \, \ \ \forall w\in J, \, \forall n\in\Z.
\eequ
\end{lemma}
\bdem If $\varphi$ comes from \eqref{mod-h},
$$
|\varphi^{(n)'}(w)| = \frac{  1-
\left|\frac{1-\alpha^n}{1+\alpha^n}\right|^2  }{ \left|1-
\frac{1-\alpha^n}{1+\alpha^n}w\right|^2 } = \frac{
\frac{4\alpha^n}{(1+\alpha^n)^2}  }{ \left|1-
\frac{1-\alpha^n}{1+\alpha^n}w\right|^2 }.
$$
When $w\in J$ and $n\in \Z$, $\frac{1-\alpha^n}{1+\alpha^n} w$
remains in the angular sector
$\{r e^{i\theta}: 0 \le r \le 1 , \ \arg \varphi(i) < \theta <\arg \varphi(-i) \}.$
Since $\disc \cap \partial D_\varphi$ is orthogonal to $\partial \disc$ (see Figure 1), it follows that
$$
\frac{2^2}{(\alpha+1)^2}=|1-\varphi(0)|^2 \le \left|1-
\frac{1-\alpha^n}{1+\alpha^n}w\right|^2 \le 2^2 .
$$
For $n\in\Z$,
$$
\alpha^{-|n|} \le \frac{4\alpha^n}{(1+\alpha^n)^2} \le 4
\alpha^{-|n|} ,
$$
so combining the above inequalities we obtain the desired result.
If $\varphi$ comes from \eqref{mod-p},
$$
|\varphi^{(n)'}(w)| = \frac{  1- \left|\frac{nt}{nt-2i}\right|^2
}{ \left|1- \frac{nt}{nt-2i}w\right|^2 } =\frac{  \frac{4}{n^2t^2
+4}  }{ \left|1- \frac{nt}{nt-2i}w\right|^2 }.
$$
When $w\in [-1,\varphi(-1))$, Figure 2 shows that
\begin{eqnarray*}
2\ge \left|\ov{w}- \frac{nt}{nt-2i}\right| &\ge |\varphi(-1)-1|-
\left|1-\frac{nt}{nt-2i}\right| \\*[1mm] &= 2 \left[\frac{1}{
(t^2+1)^{1/2} } - \frac{1}{ (n^2t^2+4)^{1/2} }\right]  \\*[1mm]
&\stackrel{\mbox{\scriptsize{if $n\neq 0$}}}{\ge} 2
\left[\frac{1}{ (t^2+1)^{1/2} } - \frac{1}{ (t^2+4)^{1/2} }\right] .
\end{eqnarray*}
The claim follows for $n\neq 0$ by inserting the above inequalities in
the expression of $|\varphi^{(n)'}(w)|$. Since $|\varphi^{(0)'}(w)|=1$, it also follows for $n=0$ by taking $c(t)>4$.
\edem

\begin{lemma}\laba{cutout}
Let $1\le p<\infty$ and $\varphi$ be a non-elliptic automorphism.
If $f_0 \ge 0$ belongs to $L^p(J)$ and\/
$\log f_0\in L^1 (J)$ with respect to the linear Lebesgue measure
$|dz|$, then for $\lambda >0$,
$$
f(z) := \sum_{n\in \Z} \chi_{ \varphi^{(n)}(J) }(z)  \, \lambda^n
f_0( \varphi^{(-n)}(z) )
$$
is in $L^p(\partial \disc)$ and\/ $\log f\in L^1 (\partial \disc)$
if and only if
\begin{enumerate}
\item[(1)] $\frac{1}{\sqrt[p]{\alpha}}< \lambda < \sqrt[p]{\alpha}\,$ in the hyperbolic case (i.e.: when $\varphi$ is given by \eqref{mod-h}),
\item[(2)] $\lambda =1$ in the parabolic case (i.e.: when $\varphi$ is given by \eqref{mod-p}).
\end{enumerate}
Furthermore, $f\circ \varphi = \lambda f$ on $\partial \disc$ in
both cases.
\end{lemma}
\bdem Clearly $f$ is defined almost everywhere on $\partial \disc$ so that $f|_J = f_0$ and $(f\circ \varphi^{(n)})(z) = \lambda^n
f(z)$ for all $n\in\Z$ (see Figures 1 and 2).
Thus,
\begin{eqnarray}
\int_{ \varphi^{(n)}(J) } f(z)^p |dz|
&= & \displaystyle  \int_{J} f(\varphi^{(n)}(w))^p |\varphi^{(n)'}(w)| |dw| \nonumber \\
&= & \displaystyle \int_{J} \lambda^{pn} f(w)^p |\varphi^{(n)'}(w)| |dw|.             \laba{ei1}
\end{eqnarray}
In the hyperbolic case we insert the inequalities \eqref{des-h} in the above expression to obtain  a constant $C(\alpha) >0$ such that
$$
\frac{1}{4} \frac{ \lambda^{pn} }{ \alpha^{|n|} } \int_{J} f(w)^p
|dw| \leq \int_{\varphi^{(n)}(J)} f(z)^p \, |dz| \leq C(\alpha)
\frac{ \lambda^{pn} }{ \alpha^{|n|} } \int_{J} f(w)^p  |dw| .
$$
Summing over $n\in\Z$,
$$
\frac{1}{4} \sum_{n\in\Z} \frac{ \lambda^{pn} }{ \alpha^{|n|} }
\int_{J} f^p \leq \int_{\partial\disc} f^p \leq C(\alpha)
\sum_{n\in\Z} \frac{ \lambda^{pn} }{ \alpha^{|n|} } \int_{J} f^p .
$$
For the convergence of the sum it is necessary and sufficient that
$\lambda^p/ \alpha <1$ and $\lambda^p \alpha >1$; that is,
$\frac{1}{\sqrt[p]{\alpha}}< \lambda < \sqrt[p]{\alpha}$. The same
idea, with $f^p$ replaced by $\log f$, gives
$$
\frac{1}{4} \sum_{n\in\Z} \frac{ 1 }{ \alpha^{|n|} } \int_{J}
(\log f + n\log\lambda) \leq \int_{\partial\disc} \log f \leq
C(\alpha)  \sum_{n\in\Z} \frac{ 1 }{ \alpha^{|n|} } \int_{J} (\log
f + n\log\lambda),
$$
which converges because $\alpha >1$.

\noi In the parabolic case we insert the inequalities
\eqref{des-p} in \eqref{ei1} to obtain a constant $c(t) >0$ such
that
$$
\frac{ \lambda^{pn} }{ n^2 t^2+4 } \int_{J} f(w)^p  |dw| \leq
\int_{\varphi^{(n)}(J)} f(z)^p \, |dz| \leq c(t)\frac{
\lambda^{pn} }{ n^2 t^2+4 } \int_{J} f(w)^p  |dw| .
$$
Summing over $n\in\Z$,
$$
\sum_{n\in\Z} \frac{ \lambda^{pn} }{ n^2 t^2+4 } \int_{J} f^p \leq
\int_{\partial\disc} f^p \leq c(t)\sum_{n\in\Z} \frac{
\lambda^{pn} }{ n^2 t^2+4 } \int_{J} f^p .
$$
The sum  converges if and only if $\lambda=1$. Finally, observing that for
$\lambda=1$ we have $\log (f\circ \varphi) =\log f$, the same
proof, with $f^p$ replaced by $\log f$, gives
$$
\sum_{n\in\Z} \frac{ 1 }{ n^2 t^2+4 } \int_{J} \log f \leq
\int_{\partial\disc} \log f \leq c(t) \sum_{n\in\Z} \frac{ 1 }{
n^2 t^2+4 } \int_{J} \log f ,
$$
which converges.
\edem

\begin{theo}\laba{caout}
Let $1\le p<\infty$ and $\varphi$ be a non-elliptic automorphism.
There is an outer function
$$
F(z)= \exp \left( \int \frac{\eiti +z}{\eiti -z}  \log f(\eiti) \,
\frac{d\theta}{2\pi} \right) \in H^p
$$
such that $F\circ \varphi = \gamma\lambda  F$ for $\lambda >0$ and
some $\gamma \in\partial\disc$ if and only if
\begin{enumerate}
\item[(1)] $\frac{1}{\sqrt[p]{\alpha}}< \lambda < \sqrt[p]{\alpha}$ in the hyperbolic case (where $\varphi$ is given by \eqref{mod-h}),
\item[(2)] $\lambda =1$ in the parabolic case,
\end{enumerate}
and $f$ is given by Lemma \ref{cutout}. Moreover, $\gamma$ depends
on $\varphi$ and $f_0$ but it is independent of $\lambda$ in both
cases.
\end{theo}
\bdem Since $|F(\xi)|= f(\xi)$ for almost every
$\xi\in\partial\disc$, we also have $|F(\varphi(\xi))|=
f(\varphi(\xi))=\lambda f(\xi)$, where the last equality comes
from the lemma. That means that the outer functions $F\circ \varphi$ and $\lambda F$ have the same modulus on
$\partial \disc$, and consequently they differ by a multiplicative constant of modulus one (see \cite[Ch.$\,$II, Thm.$\,$4.6]{Ga}).
To complete the proof, notice that
\begin{eqnarray*}
\gamma\lambda = \frac{F(\varphi(0))}{F(0)}
&= &
\exp \left( \int_{-\pi}^{\pi} \left[\frac{\eiti +\varphi(0)}{\eiti -\varphi(0)} -1\right]
 \log f(\eiti) \frac{d\theta}{2\pi} \right) \\*[1mm]
&= &
\exp \left( \int_{\partial\disc} \frac{2\varphi(0)}{\eiti-\varphi(0)}  \log f(\eiti) \frac{|d(\eiti)|}{2\pi} \right) .
\end{eqnarray*}
So, if we partition $\partial\disc = \bigcup_{n\in \Z} \varphi^{(n)}(J)$
and change variables (as in \eqref{ei1}), we obtain
\begin{eqnarray*}
\gamma &= &
\exp \left(  i\!\int_{\partial\disc} \Im\frac{2\varphi(0)}{\eiti-\varphi(0)}  \log f(\eiti) \frac{|d(\eiti)|}{2\pi} \right) \\*[1mm]
&= &    \exp \left( i 2\varphi(0) \int_J \sum_{n\in\Z}\Im\frac{ |\varphi^{(n)'}(\eiti)| }{
\varphi^{(n)}(\eiti) -\varphi(0) } \, [\log f_0(\eiti)+ n\log \lambda]  \, \frac{|d(\eiti)|}{2\pi} \right) \\*[1mm]
&= &
   \exp \left( i 2\varphi(0) \int_J \log f_0(\eiti) \sum_{n\in\Z}\Im\frac{ |\varphi^{(n)'}(\eiti)| }{
\varphi^{(n)}(\eiti) -\varphi(0) } \,\frac{|d(\eiti)|}{2\pi} \right) ,
\end{eqnarray*}
where the last equality holds in the hyperbolic case because $J$
is symmetric with respect to the real axis,
$\varphi(0)\in\mathbb{R}$ and $\varphi^{(n)}(e^{-i\theta}) = \ov{
\varphi^{(n)}(\eiti) }$, and in the parabolic case because
$\lambda=1$.
\edem

\begin{rema}\laba{rumia}
{\em The case $p=\infty$ is easier. If a bounded function
is an eigenvector corresponding to an eigenvalue $a$, then  $|a|= 1$. On the other hand, if $f_0\in L^\infty(J)$, $\log f_0\in
L^1(J)$ and $\lambda=1$ in Lemma \ref{cutout}, it is clear that
the function $f$ of the lemma is bounded and the same proof shows
that $\log f\in L^1(\partial\disc)$. Thus, Theorem \ref{caout}
holds for $p=\infty$ with $\lambda=1$ in both the hyperbolic and
parabolic case. }
\end{rema}

\subsection{Singular inner eigenvectors}
For $\mu$ a positive, finite, singular measure on $\partial \disc$, its associated singular inner function is
$$
S_\mu(z)= \exp \left( -\int \frac{\eiti +z}{\eiti -z}  \, d\mu(\eiti)  \right).
$$

The following lemma was proved by  Matache in \cite{ma} as a means to characterize
the singular inner eigenfunctions in terms of pull-back measures.
We give a different proof.

\begin{lemma}
Let $\varphi$ be any automorphism and $\nu$ be a (finite positive) singular measure on $\partial \disc$.
Then there exists a unique singular measure $\mu$ that satisfies
$|S_\nu(\varphi^{-1}(z))|=|S_\mu(z)|$ for $z\in\disc$, and it is given by
\bequ\laba{mata}
\mu(E)= \int_{\varphi^{-1}(E)} |\varphi'| d\nu
\eequ
for Borel sets $E\subset\partial\disc$.
\end{lemma}
\bdem
The existence holds because $S_\nu(\varphi^{-1}(z))$ is a singular inner function, so there must exist a singular measure $\mu$ and
$\lambda\in\partial\disc$ such that $S_\nu(\varphi^{-1}(z))= \lambda S_\mu(z)$. The uniqueness follows because two singular inner functions
with the same modulus have the same associated measure (see \cite[p.$\ $70]{Ga}). To prove the last statement consider first the case  $\nu = \delta_\xi$,
the Dirac measure concentrated at a point $\xi\in\partial\disc$.
If $\mu_\xi$ is the measure that satisfies $|S_{\delta_\xi}(\varphi^{-1}(z))|=|S_{\mu_\xi}(z)|$, the first
function extends continuously to $\partial \disc \setminus \{ \varphi(\xi) \}$ as the constant 1, and thus the same holds for $|S_{\mu_\xi}(z)|$, which means
that ${\mu_\xi} = c \delta_{\varphi(\xi)}$ for some constant $c>0$.
Moreover,
$$
c = {\mu_\xi}(\disc)= -\log |S_{\mu_\xi}(0)|= -\log
|S_{\delta_\xi}(\varphi^{-1}(0))| =
\frac{1-|\varphi^{-1}(0)|^2}{|\xi-\varphi^{-1}(0)|^2} =
|\varphi'(\xi)| .
$$
So, $\mu_\xi=|\varphi'(\xi)| \delta_{\varphi(\xi)}$, which is the measure defined by \eqref{mata} for $\nu = \delta_\xi$.
For an arbitrary singular measure $\nu$, write $\nu = \int \delta_\xi \, d\nu(\xi)$ and consider the singular measure
%%  $\mu:= \int |\varphi'(\xi)| \delta_{\varphi(\xi)} d\nu (\xi)$,
$\mu:= \int \mu_\xi d\nu (\xi)$,
where the integrals converge weak-$\ast$ in the space of finite Borel measures.
Since $\mu_\xi$ is the measure that satisfies \eqref{mata} for $\delta_\xi$, then $\mu$ is the measure that satisfies \eqref{mata} for $\nu$.
Moreover, since the map $\sigma \mapsto P_z(\sigma)= -\log |S_\sigma(z)|$ (the Poisson integral of $\sigma$) is linear,
$$
 -\log |S_\nu(\varphi^{-1}(z))| =
P_{\varphi^{-1}(z)} (\nu) = \int P_{\varphi^{-1}(z)}(\delta_\xi) d\nu(\xi) = \int P_z (\mu_\xi) d\nu(\xi) = P_z(\mu) = -\log |S_\mu(z)|,
$$
where we can take the integral of measures outside of the Poisson integral because the kernel of $P_w$ is continuous on $\partial \disc$ for every $w\in \disc$, and
the equality in the middle is proved above.
\edem

In order to present a statement that is as clear as possible, we
allow the possibility of $\nu\equiv 0$ as a singular measure in
the next corollary, and we interpret this to mean that $S_\nu \equiv 1$.
\begin{coro}\laba{fixy}
Let $\varphi$ be a non-elliptic automorphism and $\nu$ be a
non-negative measure supported on the fixed points of $\varphi$.
Then\/ $|S_\nu(\varphi^{-1}(z))|=|S_\nu(z)|$ for $z\in\disc\,$ if and only if
\begin{enumerate}
\item[ (1)] $\nu \equiv 0$ in the hyperbolic case,
\item[ (2)] $\nu = a \delta_1$ in the parabolic case, where $\delta_1$
is the Dirac measure at $1$ and\/  $a\geq 0$.
\end{enumerate}
\end{coro}
\bdem By \eqref{mata}, $|S_\nu(\varphi^{-1}(z))|=|S_\nu(z)|$ if
and only if $\nu(\{\xi\})=  |\varphi'(\xi)| \nu(\{\xi\})$ for any
fixed point $\xi$ of $\varphi$. So, $\xi\in \{ -1,1\}$ when
$\varphi$ is hyperbolic, and since $|\varphi'(1)| = \alpha^{-1}$
and $|\varphi'(-1)| = \alpha$, where $\alpha\neq 1$, then
$\nu(\{\xi\})=0$. If $\varphi$ is parabolic, we have that $\xi=1$
and $|\varphi'(1)|=1$, which means that $\nu(\{1\})$ can be any
non-negative number.
\edem

Now we are in position to state the result concerning singular
inner functions.

\begin{theo}\laba{casin}
Let $\varphi$ be a non-elliptic automorphism and $J\subset
\partial \disc$ be the set associated with $\varphi$ by the
paragraph that precedes Lemma \ref{deriv}. If $\nu_0$ is a
singular measure on $\partial\disc$ with mass concentrated in $J$, write
\bequ\laba{rut}
d\nu := \sum_{n\in \Z} d\nu_n,    \ \mbox{ where }\ \nu_n(G) = \int_{\varphi^{(-n)}(G)}
|\varphi^{(n)'}(\eiti)| \, d\nu_0(\eiti) .
\eequ
Then $S_\nu$ is a singular inner function that satisfies $S_\nu(\varphi(z))= \gamma
S_\nu(z)$ for some $\gamma\in\partial\disc$. If $\varphi$ is
parabolic, the same holds for $S_{\nu + a\delta_1}$ for any $a\ge
0$. Conversely, every singular inner eigenvector $S_\mu$ has this
form, with $\nu_0= \chi_{J} \mu$ and, if $\varphi$ is
parabolic, $a=\mu(\{1\})$.
\end{theo}
\bdem First we prove that $\nu$ is a finite measure. Observe that
the whole mass of $\nu_n$ is concentrated in $\varphi^{(n)}(J)$
and that these sets are pairwise disjoint. Hence,
$$
\nu(\partial\disc)= \sum_{n\in\Z}\nu_n(\partial\disc) =
\int_{\partial\disc} \sum_{n\in\Z}|\varphi^{(n)'}(\eiti)| \,
d\nu_0(\eiti) = \int_{J} \sum_{n\in\Z}|\varphi^{(n)'}(\eiti)| \,
d\nu_0(\eiti) .
$$
It follows from Lemma \ref{deriv} that the above quantity is
finite.

By \eqref{mata}, $|S_{\nu_0}(\varphi^{(-n)}(z))|=|S_{\nu_n}(z)|$,
and consequently,

$$
|S_{\nu_{n-1}}(z)|=|S_{\nu_0}(\varphi^{(-[n-1])}(z))|=
|S_{\nu_0}(\varphi^{(-n)}(\varphi(z)))|=|S_{\nu_n}(\varphi(z))| .
$$
Therefore

$$
\log |S_\nu(z)| = \sum_{n\in\Z} \log |S_{\nu_{n-1}}(z)|
=\sum_{n\in\Z} \log |S_{\nu_n}(\varphi(z))| = \log
|S_\nu(\varphi(z))|,
$$
which means that $S_\nu(\varphi(z))/S_\nu(z)$ is some constant of
modulus $1$. Corollary \ref{fixy} says that we can add to $\nu$ an
atom at $1$ when $\varphi$ is parabolic.

To prove the converse, suppose that $S_\mu$ is a singular inner
eigenvector and write $\nu_n = \chi_{\varphi^{(n)}(J)} \mu$. If we
fix $n\in\Z$, for any Borel set $G\subset\varphi^{(n)}(J)$, the
equality $|S_\mu(\varphi^{-(n)}(z))|=|S_\mu(z)|$ together with
\eqref{mata} yields
$$
\nu_n(G) =\mu(G) = \int_{\varphi^{(-n)}(G)}
|\varphi^{(n)'}(\eiti)| \, d\mu(\eiti) = \int_{\varphi^{(-n)}(G)}
|\varphi^{(n)'}(\eiti)| \, d\nu_0(\eiti) .
$$
Thus, $\mu = \nu + \mu_\infty$, where $\nu$ is given by
\eqref{rut}, with $\nu_0 = \chi_J \mu$ and $\mu_\infty$ supported
on the fixed points of $\varphi$. But since $S_\nu$ and $S_\mu$
are both eigenvectors, the same holds for $S_{\mu_\infty}$. The
theorem now follows from Corollary \ref{fixy}.
\edem

We can summarize this discussion by saying that if $1\le p<\infty$, and $h=FBS\in H^p$, where $F$ is outer,
$B$ is a Blaschke product and $S$ is a singular inner function, then $h$
is an eigenvector of $C_\varphi$ if and only if $F$, $B$ and $S$
are respectively given by Theorems \ref{caout}, \ref{cabla} and
\ref{casin}. The same holds for $p=\infty$, with $F$ given by
Remark \ref{rumia}. \\

\noindent {\bf Acknowledgements:}
We thank the referee for carefully reading the paper and for suggestions that have helped us to improve
the exposition of the paper.
The first author is partially supported by
the grant MTM2010-16679 and the Gobierno de Arag\'on research group
\emph{An\'alisis Matem\'atico y Aplicaciones}, ref. DGA E-64.
The third author is partially supported by the Ram\'on y Cajal program
and the grants 2009SGR00420 and MTM2008-00145 (Spain), and BID-PICT2009-0082 (Argentina).
Part of this work was done while the second author was
visiting the University of Zaragoza with a grant of the research
institute IUMA.

\bibliographystyle{amsplain}

\begin{thebibliography}{99}
\bibitem{Ah} L. V. Ahlfors, Complex Analysis, McGraw-Hill, New York, 1979.

\bibitem{BSh} L. Brown and A. L. Shields, \emph{Cyclic vector in the Dirichlet space},
Trans. Amer. Math. Soc. 285 (1984), n.$\,$1, 269--304.

%%  \bibitem{Ca} S. R. Caradus, \emph{Universal operators and invariant subspaces}, Proc. Amer. Math. Soc. 23
%%  (1969) 526--527.

\bibitem{Ci-Ro} J. A. Cima and W. T. Ross, The Backward Shift on the Hardy Space, American Mathematical Society,
Providence, 2000.

\bibitem{Cowen} C. C. Cowen, \emph{Iteration and the solution of functional equations for functions analytic in
the unit disk}, Trans. Amer. Math. Soc. 265 (1981), 69--95.

\bibitem{Co-Mac} C. C. Cowen and B. D. MacCluer, Composition Operators on Spaces of Analytic Functions,
CRC Press, Boca Raton, 1995.


\bibitem{Du} P. Duren, Theory of $H^{p}$ spaces, Academic
Press, New York (1970).

\bibitem{Ga} J. B. Garnett, Bounded analytic functions,  Revised first edition.
Graduate Texts in Mathematics, 236. Springer, New York, 2007.

\bibitem{GG} E. A. Gallardo-Guti\'errez and P. Gorkin, \emph{Cyclic Blaschke products for composition
operators}, Rev. Mat. Iberoamericana 25 (2009), n.$\,$2, 447--470.

\bibitem{Hedenmalm} H. Hedenmalm, \emph{Thin interpolating sequences and three
algebras of bounded functions},  Proc. Amer. Math. Soc.  99 (1987),
n.$\,$3, 489--495.

\bibitem{ma} V. Matache, \emph{The eigenfunctions of a certain composition
operator}, Contemp. Math. 213 (1998), 121--136.

\bibitem{Mortini} R. Mortini,
\emph{Cyclic subspaces and eigenvectors of the hyperbolic
composition operator}, Travaux math\'{e}matiques, Fasc. VII, 69--79,
S\'{e}m. Math. Luxembourg, Centre Univ. Luxembourg, Luxembourg, 1995.

\bibitem{NRW1} E. A. Nordgren, P. Rosenthal and F. S. Wintrobe,
\emph{Composition operators and the invariant subspace problem},
C. R. Mat. Rep. Acad. Sci. Canada,  6, (1984), 279--282.

\bibitem{NRW} E. A. Nordgren, P. Rosenthal and F. S. Wintrobe
\emph{Invertible composition operators on $H^p$}, J. Functional
Analysis,  73 (1987), 324--344.

%%  \bibitem{RR} H. Radjavi and P. Rosenthal, Invariant subspaces,
%%  Second edition. Dover Publications, Inc., Mineola, NY, 2003.

\bibitem{Ro} G. C. Rota, \emph{On models for linear operators}, Comm. Pure Appl.
Math. 13 (1960) 469--472.

\bibitem{Su-Wo} C. Sundberg and T. H. Wolff, \emph{Interpolating sequences for $QA_B$}, Trans. Amer. Math. Soc. 276 (1983), 551–-581.
\end{thebibliography}

\end{document}